\documentclass[12pt,a4paper,twoside]{amsart}
\usepackage[english]{babel}
\usepackage{amsmath}
\usepackage{amssymb}
\usepackage[latin1]{inputenc}  
\theoremstyle{plain}
\newtheorem{theorem}{Proposition}[section]
\newtheorem{defi}[theorem]{Definition}
\newtheorem{thm}[theorem]{Theorem}

\newtheorem{corol}[theorem]{Corollary}
\renewcommand{\H}{\mathcal H}

\newcommand{\A}{\mathcal A}
\newcommand{\D}{\mathcal D}
\newcommand{\F}{\mathcal F}
\newcommand{\B}{\mathcal B}
\newcommand{\la}{\big\langle}
\newcommand{\ra}{\big\rangle}
\renewcommand{\L}{\mathcal L^+(\D)}
\begin{document}

\title[]{Jordan Maps and Triple Maps on Algebras of Unbounded Operators}
\author[]{WERNER TIMMERMANN}

\address{Institut f\"ur Analysis\\
         Technische Universit\"at Dresden\\
         D-01062 Dresden, Germany}
\email{timmerma@math.tu-dresden.de}
\begin{abstract}

The paper contains results on the structure of Jordan maps and several kinds of triple maps on standard algebras of unbounded operators in Hilbert spaces. These results are unbounded counterparts to results on algebras of bounded operators obtained by Lu, Moln\'ar and others ([3,4,6 -- 10]).

\end{abstract}

\thanks{\hspace{-0.5cm}Keywords: algebras of unbounded operators, Jordan maps, triple maps\\
 2000 Mathematics Subject Classification: 47L60, 46K05 }
\maketitle

\section{Introduction }
In the last 15 years there appeared a lot of papers with special structural results concerning mappings between certain classes of operators on Banach or Hilbert spaces. One striking point in the proofs of several of these results is the application of theorems from abstract ring theory. Let us mention two already almost classical results applied at several places.\\
In \cite{mart69} Martindale showed that every multiplicative bijective mapping from a prime ring containing a nontrivial idempotent onto an arbitrary ring is necessarily additive, and hence it is a ring isomorphism.\\
A result of Herstein \cite{herstein} (see also \cite{palmer94}) states that a Jordan homomorphism of an arbitrary ring onto a prime ring of characteristic different from 2 and 3 is a ring homomorphism or a ring antihomomorphism.\\
It turned out that some of the results for algebras of bounded operators are also valid for algebras of unbounded operators on dense domains in Hilbert spaces (see for example \cite{ti93,ti00b, ti01a, ti03a, ti04a, ti05a, ti05b}). The proofs use some of the abstract ring theoretical results as well as appropriate adaptions of functional analytical ideas used formerly in the proofs for algebras of bounded operators.\\
The aim of the present paper is to discuss Jordan maps and several kinds of triple maps on algebras of unbounded operators. The considerations were stimulated by the papers of Moln\'ar and others (see for example \cite{molnar97a,molnar98d,molnar98g, molnar02f}) . L. Moln\'ar turned my attention also to the papers of Lu (inspired in their turn by Moln\'ar's results).\\
Let us fix our notions and notation.\\
For a Banach space $\mathcal X$ let us denote by $\mathcal B(\mathcal X)$ Banach algebra of all bounded operators on  $\mathcal X$ and $\mathcal F(\mathcal X) \subset \mathcal B(\mathcal X) $ the ideal of all finite rank operators. An algebra $\mathcal A \subset \mathcal B(\mathcal X)$ is said to be a standard operator algebra if it contains $\mathcal F(\mathcal X)$. The same notions apply if $\mathcal X$ is a Hilbert space $\mathcal H$. Then $\B(\H)$, $\F(\H)$ are *-algebras and standard (*-)algebras can be defined alternatively to be subalgebras or *-subalgebras of $\B(\H)$ containing $\F(\H)$. \\
Concerning algebras of unbounded operators we use the following notions 
 (a standard reference for algebras of unbounded operators is \cite{schm90}).\\
Let $\mathcal D$ be a dense linear subspace in a Hilbert space $\mathcal H$ with scalar
product $\la ~,~ \ra $ (which is supposed to be conjugate linear in the first and linear in the second component). The set of linear operators $\L = \{ A: A \D \subset \D, A^* \D \subset \D\}$ is a $\ast$-algebra with respect to the natural operations and the involution $A \rightarrow A^+ = A^*|\D$.
The graph topology $t$ on $\D$ induced by $\L$ is generated by the directed family of seminorms $\phi \rightarrow ||\phi ||_A = || A \phi ||,~ \forall ~  A \in \L, ~  \phi \in \D$. $\D$ is called an (F)-domain, if $(\D, t)$ is an (F)-space.  
A standard  operator (*-)algebra is a (*-)subalgebra $\mathcal A(\D) \subset \L$ containing the ideal $\mathcal F(\D) \subset \L$ of all finite rank operators on $\D$.\\
Recall that every rank-one operator $F \in \F(\H)$ has the form $F = \phi \otimes \psi$, i.e. $F \chi = \la \psi, \chi\ra \phi$, where $\phi,\psi,\chi \in \D$. Every standard operator algebra (in the bounded as well as in the unbounded case) is prime, that is $A\A B = \{0\}$ implies $A = 0$ or $B = 0$. \\
Remark that a standard operator algebra $\A(\D)$ is dense in $\L$ with respect to the strong operator topology. More exactly, for every $A \in \L$ there is a net $(T_\alpha) \subset \A(\D)$ such that $T_\alpha\phi \rightarrow A\phi$ für all $\phi \in \D$. In particular, every standard operator algebra contains a (norm-bounded) approximation of the unit with respect to the strong operator topology (take the net of all ortho-projections from $\F(\D)$). \\
We allow the following slight abuse of language and notation. If $A \in \L$ is a bounded operator, then we will in general not distinguish $A$ and its closure $\bar{A} \in \B(\H)$. So we speak of a unitary operator $U \in \L$ meaning that $\bar{U}$ is a unitary operator in $\B(\H)$.\\
Some words about the content of the paper.\\
In section 2 we discuss the structure of Jordan maps between standard operator algebras. We use the (mostly algebraic) results of Lu \cite{lu03b} to get essentially the same structural results of these mappings as Moln\'ar \cite{molnar02f}.\\
In section 3 we describe the structure of several triple maps on standard operator algebras $\A(\D)$. The corresponding results in the bounded case are given in [5,10,11]. \\
The specific notions and notation will be given in the corresponding sections. 
\section{Jordan maps}
Let $\mathcal R, \mathcal R'$ be rings. A mapping $\Phi\colon \mathcal R \rightarrow \mathcal R'$ is called a Jordan homomorphism if $\Phi$ is additive and satisfies
\begin{equation}
\Phi(A^2) = \Phi(A)^2 \quad (A \in \mathcal R).
\end{equation}
If $\mathcal R'$ is 2-torsion free (i.e. $2A = 0$ implies $A = 0$) then (1) is equivalent to
\begin{equation}
\Phi(AB + BA) = \Phi(A) \Phi(B) + \Phi(B) \Phi(A) \quad (A,B \in \mathcal R).
\end{equation}
Typical examples of Jordan homomorphisms are ring homomorphisms and ring antihomomorphisms (i.e. additive mappings satifying $\Phi(AB) = \Phi(B) \Phi(A)$).\\
In algebras one considers besides (2) the equation
\begin{equation}
\Phi\big(\frac{1}{2}(AB + BA)\big) = \frac{1}{2} \big[\Phi(A) \Phi(B) + \Phi(B) \Phi(A)\big].
\end{equation}
Clearly, if $\Phi$ is additive then (2) and (3) are also equivalent.\\
The expression $A\circ B := AB + BA$ is often called the Jordan product of $A,B$. So, a Jordan homomorphism is an additive mapping that preserves the Jordan product. The interest in studying such maps was essentially stimulated by the mathematical foundation of quantum mechanics. There the observables are represented by selfadjoint operators from $\B(\H)$. The subspace of all selfadjoint operators from $\B(\H)$ is not closed under the usual product but under the Jordan product.\\In \cite{lu03b} Lu introduced a more general notion. Let $\A , \A'$ be algebras over the field $\mathbb Q$ of rational numbers. A bijective $\Phi\colon \A \rightarrow \A'$ such that
\begin{equation}
\Phi(k(AB + BA)) = k\big[\Phi(A)\Phi(B) + \Phi(B)\Phi(A)\big]
\end{equation}
is called a k-Jordan map. Emphasize that no additivity is assumed!\\
Again, additivity of $\Phi$ implies that (1) and (4) are equivalent.\\
In \cite{molnar02f} Moln\'ar clarified the structure of $\frac{1}{2}$-Jordan maps between standard operator algebras on Banach spaces (Theorem 1), $\frac{1}{2}$- Jordan *-maps between standard operator *- algebras on Hilbert spaces (*-map means $\Phi(A^*) = \Phi(A)^*$) (Corollary 2) and 1-Jordan *-maps between standard operator *-algebras on Hilbert spaces (Theorem 3). He left open the case of 1-Jordan maps between standard operator algebras on Banach spaces. Let us summarize  his Corollary 2 and Theorem 3.\\

{\bf Theorem A}(\cite{molnar02f}) \\
\emph{ Let $\H , \mathcal K$ be Hilbert spaces, dim $ \H  > 1$, and let $\A \subset \B(\H), \B \subset \B(\mathcal K)$ be standard operator *-algebras. Suppose that $\Phi\colon \A \rightarrow \B$ is a bijective 1- Jordan (or $\frac{1}{2}$ -Jordan) *-map. Then we have the following four possibilities:\\
i) there exists an either unitary or antiunitary operator $U \colon \H \rightarrow \mathcal K$ such that
$$\Phi(A) = UAU^* \quad (A \in \A);$$
ii) there exists an either unitary or antiunitary operator $U\colon \H \rightarrow \mathcal K$ such that
$$\Phi(A) = UA^*U^* \quad (A \in \A). $$}

The result of Lu \cite{lu03b}  is quite general and solved also the open problem in \cite{molnar02f} mentioned above. Lu's abstract algebraic result reads as follows.\\

{\bf Theorem B}( \cite{lu03b}, Theorem 1.1)\\
\emph{Let $\A, \A'$ be algebras over $\mathbb Q$. Suppose that $\A$ contains an idempotent $e_1$ which satisfies:\\
i) $ e_ixe_j\A e_l = 0$ or $e_l\A e_ixe_j = 0$ implies $e_ixe_j = 0 ~(1 \leq i,j,l \leq 2)$, where $e_2 = 1 - e_1$ ($\A$ need not have an identity element, so $e_2$ need not belong to $\A$).\\
ii) if $e_2te_2xe_2 + e_2xe_ste_2 = 0$ for every $t \in \A$, then $e_2te_2 = 0$.\\
Let $k$ be a fixed nonzero rational number. Suppose that $\Phi\colon \A \rightarrow \A'$ is a k-Jordan map satisfying\\
iii) $\Phi(e_1ae_2 + e_sbe_1) = \Phi(e_1ae_2) + \Phi(e_2be_1)$ for all $a,b \in \A$.\\}
\emph{Then $\Phi$ is additive.}\\

The proof of this theorem is rather long and a clever application of some decomposition ideas of Martindale \cite{mart69}. Theorem B is applied to prove\\

 {\bf Theorem C} ( \cite{lu03b}, Theorem 1.6)\\
\emph{Let $\mathcal X$ be a Banach space, dim $\mathcal X > 1$, and let $\A \subset \B(\mathcal X)$ be a standard operator algebra. Then every k-Jordan map $\Phi$ from $\A$ onto an arbitrary algebra $\A'$ over $\mathbb Q$ is additive. Moreover, $\Phi$ is either a ring homomorphism or a ring anti-homomorphism.}\\

The proof of Theorem C is purely algebraic except one step, namely to see that a standard operator algebra satisfies condition ii) in Theorem B. But this follows either directly from algebraic properties of standard operator algebras or from the existence of an approximation of the unit with respect to the strong operator topology. Moreover it is not used that $\mathcal X$ is a Banach space. \\
So, this proof applies also to standard operator algebras  $\A \subset \L$. 
The crucial point for getting the counterpart of Theorem A is the structure of ring isomorphisms or ring antiisomorphisms of standard operator algebras on $\D$. But this follows essentially from \cite{ti03a} Theorem 3.1. This theorem as it stands there is applicable only for ring isomorphisms, but a careful inspection of the proof shows that it applies also to antiisomorphisms. For sake of completeness let us formulate the results on the structure of ring isomorphisms and ring antiisomorphisms.\\

{\bf Corollary D (to \cite{ti03a} Theorem 3.1.)}
\emph{Let $\A, \B \subset \L$ be  standard operator algebras on an $(F)$-domain $\D$.} \\
\emph{i) Every ring isomorphism $\Phi: \A \to \B$ has the form
$$\Phi(A) = TAT^{-1}$$
with a bijective either linear or conjugate linear operator $T: \D \to \D$.\\
ii) Every ring antiisomorphism $\Phi: \A \to \B$ has the form
 $$\Phi(A) = TA^+T^{-1}$$}
\emph{with a bijective either linear or conjugate linear operator $T: \D \to \D$.\\
If the operator $T$ is linear, then $T \in \L$.}
\\

Now let us formulate the unbounded version of Theorem C.

\begin{thm}
Let $\D \subset \H$ be a dense linear subspace such that $(\D,t)$ is an (F)-space and let $\A, \B \subset \L$ be standard operator algebras on $\D$. Suppose that $\Phi\colon \A \to \B$ is a k-Jordan map, i.e. $\Phi$ is bijective and there is a $k \in \mathbb Q$ such that $\Phi(k(AB+BA)) = k \left[\Phi(A)\Phi(B) + \Phi(B)\Phi(A)\right]$.\\
Then $\Phi$ is a ring isomorphism or a ring antiisomorphism. So we have the following four possibilities:\\
i) there exists a bijective either linear or conjugate linear operator $T\colon \D \to \D$ such that
$$\Phi(A) = TAT^{-1} \quad (A \in \A),$$
or\\
ii) there exists a bijective either linear or conjugate linear operator $T\colon\D \to \D$ such that
$$ \Phi(A) = TA^+T^{-1} \quad (A \in \A).$$
If $T$ is linear, then it belongs to $\L$.

\end{thm}
\begin{proof}{}
From Theorem B and the proof of Theorem C combined with the density of $\A$ in $\L$ with respect to the strong operator topology it follows that $\Phi$ is additive.  Because of $k \in \mathbb Q$ we obtain that $\Phi$ is a bijective Jordan homomorphism. By the above mentioned result of Herstein \cite{herstein} $\Phi$ is a ring isomorphism or a ring antiisomorphism. Now  Corollary D applies. 
\end{proof}

\begin{corol}
Let $\D \subset \H$ be a dense linear subspace such that $(\D,t)$ is an (F)-space and let $\A, \B \subset \L$ be standard *-operator algebras on $\D$. Suppose that $\Phi\colon \A \to \B$ is a k-Jordan map and $\Phi(A^+) = \Phi(A)^+$ for all $A \in \A$. Then we have the same possibilities as in Theorem 2.1 but $T$ can be taken to be an either unitary or antiunitary operator. If $T$ is a unitary operator, then it belongs to $\L$.
\end{corol}
\begin{proof}{}
This follows from \cite{ti93} or as in \cite{molnar02f}, proof of Corollary 2.
\end{proof}

Remark: Similar results are valid if we consider maps $\Phi$ between standard operator  algebras $\A, \B$ on different domains $\D_1, \D_2$.\\

\section{Triple maps}
 There are many possibilities to define triple maps. We refer for example to the papers \cite{molnar97a,molnar97b, molnar98d}.\\
The following notions were introduced by Lu \cite{lu03a}. \\
\begin{defi}
Let $\A, \B$ be algebras (over $\mathbb Q$) and $k \in \mathbb Q$. A bijective mapping $\Phi\colon A \to B$ is called\\
 \emph{k-Jordan semi-triple map} if
$$\Phi(kABA) = k \Phi(A) \Phi(B) \Phi(A) \quad (A,B \in \A), $$
\emph{k-Jordan triple map} if
$$\Phi(k(ABC + CBA)) = k\left[\Phi(A)\Phi(B)\Phi(C) + \Phi(C)\Phi(B)\Phi(A)\right] \quad (A,B,C \in \A).$$ 
\end{defi}
Clearly, a k-Jordan triple map is a 2k-Jordan semi-triple map.\\
The main result in Lu \cite{lu03a} reads as follows.\\

{\bf Theorem E}(\cite{lu03a} , Theorem)\\
\emph{Let $\mathcal X$ be a Banach space of dimension $> 1$ and suppose that $\A$ is a standard operator algebra on $\mathcal X$. Let $\B$ be an algebra over $\mathbb Q$ and $k \in \mathbb Q$ be non-zero. Let $\Phi\colon \A \to \B$ be a bijection satisfying
$$\Phi(kABA) = k \Phi(A)\Phi(B)\Phi(A) \quad (A,B \in \A).$$
Then $\Phi$ is additive. Moreover $\Phi$ is a ring isomorphism or the negative of a ring isomorphism or a ring antiisomorphism or the negative of a ring antiisomorphism.}\\

{\bf Corollary}(\cite{lu03a}, Corollary)\\
\emph{Let $\mathcal X$ be a Banach space of dimension $> 1$ and suppose that $\A$ is a standard operator algebra on $\mathcal X$. Let $\B$ be an algebra over $\mathbb Q$ and $k \in \mathbb Q$ be non-zero. Let $\Phi\colon \A \to \B$ be a bijection satisfying
$$\Phi(k(ABC+CBA)) = k\left[\Phi(A)\Phi(B)\Phi(C) + \Phi(C)\Phi(B)\Phi(A)\right] \quad(A,B,C \in \A).$$
Then $\Phi$ is additive. Moreover, $\Phi$ is a ring isomorphism or the negative of a ring isomorphism or a ring antiisomorphism or the negative of a ring antiisomorphism.}\\

The proof of Theorem E is purely algebraic. It is used only that $\A$ is a standard operator algebra. Having the additivity of $\Phi$ one applies Theorem 3.3 in \cite{bresar89a} to conclude that $\Phi = \pm \Psi$ where $\Psi$ is a ring isomorphism or a ring antiisomorphism.\\
So the whole proof works also for standard operator algebras $\A \subset \L$. Together with an application of Corollary D we obtain the following theorem.
\begin{theorem}
Let $\D \subset \H$ be a dense linear subspace such that $(\D,t)$ is an (F)-space and let $\A, \B \subset \L$ be standard *-operator algebras on $\D$. Suppose that $\Phi\colon \A \to \B$ is a bijective mapping satisfying
$$\Phi(k(ABA) = k\Phi(A)\Phi(B)\Phi(A)  \qquad(A,B \in \A).$$
 Then $\Phi = \pm \Psi$, where $\Psi$ is a ring isomorphism or a ring antiisomorphism. Therefore $\Phi$ has one of the following  forms:\\
i)  there exists a bijective either linear or conjugate linear operator $T\colon \D \to \D$ and $c \in \{-1, 1\}$ such that
$$\Phi(A) = cTAT^{-1} \qquad (A \in \A),$$
or\\
ii) there exists a bijective either linear or conjugate linear operator $T\colon\D \to \D$ and $c \in \{-1, 1\}$ such that
$$ \Phi(A) = cTA^+T^{-1} \qquad (A \in \A).$$
If $T$ is linear it belongs to $\L$.
\end{theorem}

\begin{corol}
Let $\D \subset \H$ be a dense linear subspace such that $(\D,t)$ is an (F)-space and let $\A, \B \subset \L$ be standard *-operator algebras on $\D$. Suppose that $\Phi\colon \A \to \B$ is a bijective mapping satisfying
$$\Phi(k(ABC+CBA)) = k\left[\Phi(A)\Phi(B)\Phi(C) + \Phi(C)\Phi(B)\Phi(A)\right] \quad(A,B,C \in \A).$$
 Then the assertion of Theorem 3.2 is valid. If additionally $\Phi(A^+) = \Phi(A)^+ ~(A \in \A)$ then  $T$ can be taken to be either unitary or antiunitary. If  $T$ is linear it belongs to $\L$.
\end{corol}

Remark that the results of Lu imply the results of Moln\'ar \cite{molnar02f}. The proofs given by Moln\'ar were mainly functional analytical.
\begin{corol}
Let $\D \subset \H$ be a dense linear subspace such that $(\D,t)$ is an (F)-space and let $\A, \B \subset \L$ be standard *-operator algebras on $\D$. Suppose that $\Phi\colon \A \to \B$ is a bijective mapping satisfying
\begin{equation}
\Phi(ABC) = \Phi(A)\Phi(B)\Phi(C) \qquad (A,B,C \in \A).
\end{equation}
Then $\Phi$ has one of the following forms. There exists a bijective either linear or conjugate linear operator $T\colon \D \to \D$ and $c \in \{-1,1\}$ such that
$$\Phi(A) = cTAT^{-1} \qquad (A \in \A).$$
If $T$ is linear it belongs to $\L$.
\end{corol}
\begin{proof}{}
Equation (5) means that $\Phi$ is an 1-Jordan semi-triple map. So Theorem 3.2 applies. We show that the possibilities ii) can not occur. Suppose at the contrary that
$$\Phi(A) = cTA^+T^{-1}, \quad c \in \{-1,1\}.$$
Then on the one hand $\Phi(ABC) = cT(ABC)^+T^{-1} = cTC^+B^+A^+T^{-1}$. On the other hand due to (5)  $\Phi(ABC) = c^3 TA^+B^+C^+T^{-1}$. Therefore
$C^+B^+A^+ = A^+B^+C^+$ for all $A,B,C \in \A$ and in particular $FGH = HGF$ for all $F,G,H \in \F(\D)$ which is obviously a contradiction.
\end{proof}

Now we consider several triple maps which contain also adjoint operators. This was done  in \cite{molnar97a,molnar97b}
 and in a very general form in \cite{molnar98d}. In the first two papers the assumptions concerning $\Phi$ include now  linearity. We give two unbounded versions of such results. Because of the linearity of the maps $\Phi$ no assumptions on $\D$ are necessary.\\
We start with a theorem which is the unbounded counterpart of Proposition 1 in \cite{molnar97a}.
\begin{theorem} 
Let $\A$ be a standard operator *-algebra on $\D$. Suppose that $\Phi: \A \to \A$ is a bijective linear mapping satisfying
\begin{equation}
\Phi(AB^+C) = \Phi(A)\Phi(B)^+\Phi(C) \qquad (A,B,C \in \A)
\end{equation}
Then there are unitary operators $U,V \in \L$ such that
$$\Phi(A) = UAV \qquad (A \in \A)$$
\end{theorem}
\begin{proof}{}
We apply Theorem 3.1 in \cite{ti00b}. If we define the triple product $(A,B,C) \to AB^+C$, then $\A$ equipped with this product is a ternary ring of unbounded operators and equation (6) means that the pair $(\Phi, \Phi)$ is an automorphism pair. By Theorem 3.1 of \cite{ti00b} there are bijective linear operators $U,V \in \L$ such that
$$\Phi(A) = UAV, \quad \text{ and } \quad \Phi(A) = (U^+)^{-1} A (V^+)^{-1} \qquad (A \in \A).$$
This implies that $U,V$ are unitary operators. 
\end{proof}

The next Theorem is the unbounded version of Theorem 3 in \cite{molnar97b}.
\begin{theorem}
Let $\A, \B$ be standard operator *-algebras on $\D$. Suppose that $\Phi\colon \A \to \B$ is a bijective linear mapping such that
$$\Phi(AB^+A) = \Phi(A)\Phi(B)^+\Phi(A)\qquad (A,B \in \A).$$
Then there are either unitary operators $U,V \in \L$ such that
$$\Phi(A) = UAV \qquad (A \in \A)$$
or there are antiunitary operators $U',V'$ (mapping $\D$ onto $\D$) such that
$$\Phi(A) = U'A^+V' \qquad (A \in \A).$$
\end{theorem}
\begin{proof}{}
Because the proof is almost the same as in \cite{molnar97b} we indicate only the main steps. First, $\Phi$ presereves rank-one operators in both directions, i.e. $F \in \F$ has rank one if and only if $\Phi(F) \in \F$ has rank one.This follows because $F ( \not= 0)$ has rank one if and only if $FB^+F$ is a scalar multiple of $F$ for all $B \in \A(\D)$. By  Proposition 2.2 in \cite{ti93}  either there are bijective linear operators $S,T: \D \to \D$ such that
\begin{equation}
\Phi(\phi  \otimes \psi) = S \phi \otimes T \psi, \qquad (\phi, \psi \in  \D).
\end{equation}
or there are bijective conjugate linear operators $S', T': \D \to \D$ such that
\begin{equation}
\Phi(\phi \otimes \psi) = S'\psi \otimes T'\phi, \qquad (\phi, \psi \in \D).
\end{equation}
Suppose that $\Phi$ is of form (7) (the other case is treated similarly). Using equation $\Phi(AB^+A) = \Phi(A)\Phi(B)^+\Phi(A)$ for $A = \chi \otimes \phi, B = \rho \otimes \psi $ we get
\begin{equation}
<\phi, \psi><\rho, \chi> S\chi \otimes T\phi = <T\phi, T\psi><S\rho, S\chi>  S\chi \otimes T\phi.
\end{equation}
This leads finally to 
$$<T\phi, T\psi> = \lambda <\phi, \psi>, \qquad (\phi, \psi \in \D)$$
with some constant $\lambda$. So, $T$ is a multiple of an isometry. The same is true for $S$. Equation (9) shows that one may suppose that these multiples are equal to one. This and the bijectivity of $S,T$ imply that $U := S, V := T^+$ are unitary operators from $\L$ and
\begin{equation*}
\Phi(F) = UFV, \qquad (F \in \F(\D)).
\end{equation*}
To get $\Phi(A) = UAV$ for all $A \in \A(\D)$ is almost standard. Namely, let $A \in \A(\D), F \in \F(\D)$ be arbitrary. Then
$$UFA^+FV = \Phi(FA^+F) = \Phi(F)\Phi(A)^+\Phi(F) = UFV\Phi(A)^+UFV.$$
So $FA^+F = FV\Phi(A)^+UF$ for all $F \in \F(\D)$. This finally gives the desired result. 
\end{proof}

Remark that Theorem 3.5. can also be obtained as a immediate corollary to Theorem 3.6 below. It is easy to see that under the assumptions of Theorem 3.5. the second possibility of the form of $\Phi$ described in Theorem 3.6. can not occur (see the proof of Corollary 3.4).\\

Finally let us mention a very general result from \cite{molnar98d}. In Theorem 5 of that paper it is proved that every bijective $\Phi$ which satisfies an  n-variable *-identity is in fact an additive triple automorphism of the form (6). 
The proof is a long chain of clever algebraic calculations (not using that the operators involved are bounded). These calculations give the additivity of $\Phi$ and end with the conclusion that $\Phi$ preserves orthogonality in both directions. This means $A^*B = AB^* = 0$ if and only if $\Phi(A)^*\Phi(B) = \Phi(A)\Phi(B)^* = 0$. Since in the same paper (Theorem 4) the structure of orthogonality preserving mappings is described the result  follows. Let us formulate the unbounded version of  Theorem 5 from \cite{molnar98d}.
\begin{theorem}
Let $\D \subset \H$ be an (F)-domain and let $\A \subset \L$ be a *-ideal of $\L$. Let $\tau_j, 1 \leq j \leq n$ denote the identity on $\A$ or the adjoint operation $A \to A^+$. Suppose that $\Phi: \A \to \A$ is a bijective function satisfying the identity
$$\Phi(\tau_1(T_1)\tau_2(T_2)\dots \tau_n(T_n)) = \tau_1(\Phi(T_1))\tau_2(\Phi(T_2)) \dots \tau_n(\Phi(T_n))$$
for all $T_1,T_1,...,T_n \in \A$. 
If there is an index $j$ such that $\tau_j$ is the adjoint operation, then $\Phi$ is an additive triple automorphismus (i.e. $\Phi(AB^+C) = \Phi(A)\Phi(B)^+\Phi(C)$).\\
\end{theorem}
\begin{proof}{}
The proof up to the conclusion that $\Phi$ is additive and preserves orthogonality in both directions is exactly the same as in \cite{molnar98d}. Then one applies Theorem 3.1. from \cite{ti05b} describing the structure of bijective additive mappings preserving orthogonality in both directions. This concludes the proof.
\end{proof}


\begin{thebibliography}{13}

\bibitem{bresar89a} M. Bre\u{s}ar, {\em Jordan mappings of semiprime rings}, J. Algebra {\bf 127} (1989), 218 -- 228.

\bibitem{herstein} I. N. Herstein, {\em Jordan derivations of prime rings}, {\bf 8} (1957), 1104 -- 1110.

\bibitem{lu03a} F. Lu, {\em Jordan triple maps}, Linear Algebra Appl., {\bf 375} (2003), 311 -- 317.

\bibitem{lu03b} F. Lu, {\em Jordan maps on associative algebras}, Comm. Algebra, {\bf 31} (2003), 2273 -- 2286.

\bibitem{mart69} W. S. Martindale III, {\em When are multiplicative mappings additive?}, Proc. Amer. Math. Soc. {\bf 21} (1969), 695 -- 698.

\bibitem{molnar97a} L. Moln\'ar and B. Zalar, {\em  Three-variable *-identities and ring homomorphisms of operator ideals}, Publ. Math. Debrecen, {\bf 50} (1997), 121 -- 233.

\bibitem{molnar97b} L. Moln\'ar and P. \u{S}emrl, {\em Order isomorphisms and triple isomorphisms of operator ideals and their reflexivity}, Arch. Math., {\bf 69} (1997), 497 -- 506.



\bibitem{molnar98d} M. Györy and L. Moln\'ar and P. \u{S}emrl, {\em Linear rank and corank preserving maps on ${B}({H})$ and an application to *-semigroup isomorphisms of operator ideals }. Linear Algebra Appl., {\bf 280} (1998), 253 -- 266.


\bibitem{molnar98g}  L. Moln\'ar, {\em Characterization of additive *-homomorphisms and {J}ordan *-homomorphisms on operator ideals}, Aequationes Math., {\bf 55} (1998), 259 -- 272.

\bibitem{molnar02f} L. Moln\'ar , {\em Jordan Maps on standard operator algebras}, in: Z. Daroczy and Zs. Pales (ed.), Functional Equations - Results and Advances, Kluwer Academic Publishers, (2002), 305 -- 320.

\bibitem{palmer94} Th. W. Palmer, {\em Banach Algebras and The General Theory of *-Algebras, Vol. 1} ,Encyclopedia od Mathematics and its Applications {\bf 49} Cambridge University Press, 1994.

\bibitem{schm90} K. Schmüdgen, {\em Unbounded Operator Algebras and Representation Theory}. Akademie Verlag 1990.

\bibitem{ti93} E. Scholz and W. Timmermann, {\em Local derivations, automorphisms and commutativity preserving maps in $L^+(D)$}. Publ. RIMS, Kyoto, {\bf 29} (1993), 977 -- 995.

\bibitem{ti00b} W. Timmermann, {\em Remarks on automorphisms and derivation pairs in ternary rings of  unbounded operators}. Arch. Math. {\bf 74} (2000), 379 -- 384.

\bibitem{ti01a} W. Timmermann, {\em Additive derivations and {J}ordan derivations on algebras of  unbounded operators}. Publ. Math. Debrecen, {\bf 58} (2001), 717 -- 731.

\bibitem{ti03a} W. Timmermann, {\em Approximate derivations and isomorphisms in algebras of unbounded operators}.  Publ. Math. (Debrecen), {\bf 63} (2003),
667 -- 676.

\bibitem{ti04a} W. Timmermann, {\em Additive mappings in algebras of unbounded operators preserving operators of rank one}, Preprint TU Dresden, MATH-AN-01-2004 (2004) and arXiv math.OA/4040545. 

\bibitem{ti05a} W. Timmermann, {\em Elementary operators on algebras of unbounded operators}, Acta Math. Hungar., {\bf 107} (2005), 155 -- 166.

\bibitem{ti05b} W. Timmermann, {\em Zero product preservers and orthogonality preservers in algebras of unbounded operators }, Preprint TU Dresden, MATH-AN-02-2005 and arXiv, math.OA/0503255, 2005


\end{thebibliography}
\end{document}